\documentclass[proceedings]{aofa}

\begin{document}
	\title{Recent results on permutations without short cycles}
	
	\author{Robertas Petuchovas}
	
	\address{Vilnius University\\ Naugarduko str. 24, LT-03225 Vilnius, Lithuania\\ e-mail: robertas.petuchovas@mif.vu.lt}
	
	\maketitle

	\bigskip
		
	\newtheorem{theorem}{Theorem}
	\newtheorem{lemma}{Lemma}
	\newtheorem{corollary}{Corollary}

	\newtheorem{proposition}{Proposition}

	\def\E{\mathbf{E}}
	\def\C{\mathbf{C}}
	\def\D{\mathbf{D}}
	\def\N{\mathbf{N}}
	\def\R{\mathbf{R}}
	\def\S{\mathbf{S}_n}
	\def\Z{\mathbf{Z}}
	\def\k{\kappa}
	\def\e{\varepsilon}
	\def\n{$n\to\infty$}
	\def\cF{\mathcal F}
	
	\def\re{{\rm e}}
	\def\rO{{\rm O}}
	\def\ro{{\rm o}}
	\def\rd{{\rm d}}

	\def\s{\smallskip}
	\def\b{\bigskip}
	
	\begin{abstract} 
		{\small \textbf{Abstract.}} The density, denoted by $\kappa(n,r)$, of permutations having no cycles of length less than $r+1$ in a symmetric group $\mathrm{S}_n$ is explored. New asymptotic formulas for $\kappa(n,r)$ are obtained using the saddle-point method when $5\leq r< n$ and $n\to\infty$. 
	\end{abstract}
	

	The probability $\kappa(n,r)$ that a permutation sampled from the symmetric group $\mathrm{S}_n$ uniformly at random has no cycles of length less than $r+1$, where  $1\leq r< n$ and $n\to\infty$, is explored.
	New asymptotic formulas valid in specified regions are obtained using the saddle-point method. One of the results is  
	applied to show that estimate of the total variation distance for permutations can be expressed only through the   function $\nu(n,r)$ which is a probability that a permutation sampled from the $\mathrm{S}_n$ uniformly at random has no cycles of length greater than $r$. 
	
	To address the problem, we need recollect the following functions. Buchstab's function $\omega(v)$ is defined as a solution to difference-differential equation
	\[
	(v\omega(v))'=w(v-1)
	\]
	for $v>2$ with the initial condition $\omega(v)=1/v$ if $1\leq v\leq 2$.  Dickman's function $\varrho(v)$ is the unique continuous solution to the equation $$v\varrho'(v) + \varrho(v-1) = 0$$ for $v > 1$ with initial condition $\varrho(v)=1$ if $0\leq v\leq 1$.
	
	The interest to the problem begins with the classical example of derangements
	\begin{eqnarray*}
	\kappa(n,1)&=&\sum_{j=0}^n\frac{\left(-1\right)^{j}}{j!}=\re^{-1}+O\left(\frac{1}{n!}\right)
	\end{eqnarray*}
	and the trivial case $\kappa(n,r)=1/n$ if $n/2\leq r< n$. There was a series of works concerning general asymptotic formulas of the probability $\kappa(n,r)$ the strongest of which are presented here as Proposition 1 and Proposition 2.
	
	\begin{proposition} For $1\leq r< n$, we have
		\[
		\kappa(n,r)={\rm e}^{-H_r+\gamma}\omega(n/r)+O\left(\frac{1}{r^2}\right).
		\]
	\end{proposition}
	\noindent See~\cite[Theorem 3]{EM-on short}.

	\begin{proposition}
		Let $u=n/r$. For $1\leq r\leq n/\log n$,
		\begin{displaymath}
		\kappa(n,r)=\re^{-H_r}+O\left(\frac{(u/\re)^{-u}}{r^2}\right). 		\end{displaymath}
		If $r\geq 3$, we can replace $\re$ by 1 in the error term. 
	\end{proposition}
	
\begin{samepage}	\noindent See~\cite[Proposition 2]{AW-2015}. Together these propositions provide stronger estimates of $\kappa(n,r)$ than those in \cite{AT-AP92}, \cite{BeMaPaRi-JComTh04}, \cite{Granv-EJC06}. New results are the following theorems:
	
	\begin{theorem} For $\sqrt{n\log n}\leq r<n$, we have
		\begin{displaymath}
		\kappa(n,r)={\rm e}^{-H_r+\gamma}\omega(n/r)+O\left(\frac{\varrho(n/r)}{r^2}\right).
		\end{displaymath}
	\end{theorem}
	\noindent \textit{Proof.} The result is a corollary of \mbox{Theorem 1} in~\cite{RJ:EM-RP}. It is obtained from the probability generating function using saddle-point method, the technique is elaborated in $\cite{GT-Crible}$.
	
	\begin{theorem}
		For $(\log n)^4\leq r<n$, we have
		\begin{displaymath}
		\kappa(n,r)={\rm e}^{-H_r}+O\left(\frac{\varrho(n/r)}{r}\right).
		\end{displaymath}
	\end{theorem}
	\noindent \textit{Proof.} The saddle-point method is applied to the Cauchy's integral representation of $\kappa(n,r)$, as in the proof of Theorem 1. 
	However, there are some other technical difficulties \mbox{one must to overcome.}

	\begin{theorem}
		For $5\leq r<n$, we have
		\begin{displaymath}
		\kappa(n,r)={\rm e}^{-H_r}+O\left(\frac{\nu(n,r)}{r}\right).
		\end{displaymath}
	\end{theorem} 
	\noindent \textit{Proof.} Quite the same technique to that used in the proof of Theorem 2 is employed, just a different approximation of the saddle point is taken and Corollary~5 of \cite{EMRP-ArX15} is applied.\\
	
	Theorem 1 and Theorem 2 (see also Corollary 2.3 in~\cite{AH+GT}) improve on Proposition 1 and Proposition 2. Theorem 3 is of separate interest; as we see, it can be useful in formulas where both probabilities $\kappa(n,r)$ and $\nu(n,r)$ are involved. Here is an example.
	
	Let $k_j(\sigma)$ equal the number of cycles of length $j$ in a permutation $\sigma\in\mathrm{S}_n$, $\overline{k}(\sigma)=\left(k_1(\sigma),k_2(\sigma),\ldots,k_n(\sigma)\right),$ and $\overline{Z}=(Z_1,Z_2,\ldots,Z_n)$, where $Z_j$ are Poisson random variables such that $\mathrm{E}Z_j=1/j$, $j\in\mathrm{N}.$ Thus, if $5\leq r<n$, we have (see Lemma 3.1 on p. 69 of~\cite{ABT})
	\begin{eqnarray*}
		d_{TV}(n,r)&=&\sup_{V\subseteq \mathrm{Z}^r_+}\left|\frac{\#\{ \sigma : \overline{k}(\sigma)\in V \}}{n!}-\Pr (\overline{Z}\in V) \right|\\
		&=& \frac{1}{2}\sum_{m=0}^{\infty}\nu(m,r)\left|\kappa(n-m,r)-{\rm e}^{-H_r}\right|\\
		&=&\frac{\re^{-H_r}}{2}\sum_{m=n-r}^\infty\nu(m,r)+\frac{1}{2}\nu(n,r)+O\left(\frac{1}{r}\sum_{m=0}^{n-r-1} \nu(m,r)\nu(n-m,r)\right).
	\end{eqnarray*}
	Consequently, only results on the probability \begin{math}\nu(n,r)\end{math} are needed attempting to improve on the order of notable estimate for \begin{math} d_{TV}(n,r)\end{math} in \cite{AT-AP92}.
\end{samepage}	


\begin{thebibliography}{99}
		
		
		\bibitem{ABT}  R.~Arratia, A.D.~Barbour, and S.~Tavar\'{e}, \textit{Logarithmic Combinatorial Structures: A Probabilistic Approach}, EMS Publishing House, Z\"{u}rich, 2003.
		
		\bibitem{AT-AP92} R.~Arratia and S.~Tavar\'{e}, The cycle structure of random permutations, \textit{Ann. Probab.}, 1992, {\bf 20}, 3, 1567–-1591.
		
		\bibitem{BeMaPaRi-JComTh04} E.A.~Bender, A.~Mashatan, D.~Panario, and L.B.~Richmond, Asymptotics of combinatorial structures with large smallest component,
		\textit{ J. Comb. Th.}, Ser. A, 2004, {\bf 107}, 117--125.
		
		\bibitem{Granv-EJC06} A.~Granville, Cycle lengths in a permutation are typically Poisson, \textit{Electronic J. Comb.}, 2006, {\bf 13}, \#R107.
		
		\bibitem{AH+GT} A.~Hildebrand and G.~Tenenbaum, Integers without large prime factors, \textit{J. Th\'eorie des Nombres de Bordeaux}, 1993, {\bf 5}, 411--484.
		
		\bibitem{EM-on short} E.~Manstavi\v cius, On permutations missing short cycles, \textit{Lietuvos matem. rink.}, spec. issue, 2002, \textbf{42}, 1--6.
		
		\bibitem{RJ:EM-RP}  E.~Manstavi\v cius and R. Petuchovas, Local probabilities and total variation distance for random permutations, \textit{The Ramanujan J.}, 2016, DOI 10.1007/s11139-016-9786-0.
		
		\bibitem{EMRP-ArX15}  E.~Manstavi\v cius and R. Petuchovas, \textit{Local probabilities for random permutations without long cycles}, Electron. J. Combin., 2016, \textbf{23(1)}, \#P1.58. 
		
		\bibitem{}  E.~Manstavi\v cius and R. Petuchovas, Permutations without long or short cycles, \textit{Electronic Notes in Discrete Mathematics}, 2015, {\bf 49}, 153-158.
		
		\bibitem{GT} G.~Tenenbaum, \textit{Introduction to Analytic and Probabilistic Number Theory}, Cambridge Univ. Press, 1995.
		
		\bibitem{GT-Crible} G.~Tenenbaum, Crible d'\'{E}ratosth\`{e}ne
		et mod\`{e}le de Kubilius. In: \textit{ Number
			Theory in Progress, Proc. Conf. in Honor of Andrzej Schinzel},
		Zakopane, Poland, 1997.  K.~Gy˝ory, H.~Iwaniec, J.~Urbanowicz (Eds.), Walter de Gruyter, Berlin, New York,
		1999, 1099--1129.
		
		\bibitem{AW-2015} A.~Weingartner, On the degrees of polynomial divisors over finite fields, arXiv:1507.01920. 
		
	\end{thebibliography}
\end{document}